\documentclass{article}

\usepackage[a4paper, margin=2cm]{geometry}

\usepackage{amsthm,amsfonts,amsmath,amssymb}

\usepackage{authblk}
\author[a,b]{Petr Karnakov}
\author[a,b]{Sergey Litvinov}
\author[b,a,$\ast$]{Petros Koumoutsakos}

\affil[a]{Computational Science and Engineering Laboratory, ETH Zurich, Switzerland}
\affil[b]{School of Engineering and Applied Sciences, Harvard University, USA}
\affil[$\ast$]{corresponding author: petros@seas.harvard.edu}

\usepackage{bm}
\newcommand{\mb}[1]{\bm{#1}}
\newcommand{\pd}{\partial}

\usepackage{graphicx}
\graphicspath{{../fig/}{fig/}}
\usepackage{float}
\newcommand{\key}[1]{\raisebox{0.5pt}{\protect\includegraphics[width=5mm]{key/#1}}\,}

\usepackage{hyperref}

\title{Optimizing a DIscrete Loss (ODIL)  to solve forward and inverse problems
for partial differential equations using machine learning tools}

\begin{document}

\maketitle

\begin{abstract}
  We introduce the Optimizing a Discrete Loss (ODIL)  framework for the
  numerical solution of Partial Differential Equations (PDE) using machine
  learning tools.   The framework formulates numerical methods as a minimization
  of discrete residuals that are solved using gradient descent and Newton's
  methods. We demonstrate the value of this approach on  equations that may have
  missing parameters or where no sufficient data is available to form a
  well-posed initial-value problem.
  The framework is presented for mesh based discretizations of PDEs and inherits
  their accuracy, convergence, and conservation properties. It  preserves the
  sparsity of the solutions and  is readily applicable to inverse and ill-posed
  problems.   It is applied to  PDE-constrained optimization, optical flow,
  system identification, and data assimilation using gradient descent algorithms
  including those often deployed in machine learning.  We compare ODIL with
  related approach that represents  the solution with neural networks.  We
  compare the two methodologies and demonstrate advantages of ODIL that include
  significantly higher convergence rates and several orders of magnitude lower
  computational  cost.  We evaluate the method on various linear and nonlinear
  partial differential equations including the Navier-Stokes equations  for flow
  reconstruction problems.
\end{abstract}

\section{Main}

Classical numerical solutions for partial differential equations (PDE) rely on
their discretization, for example by finite difference~\cite{leveque2007finite}
or finite element techniques~\cite{zienkiewicz2005finite}, given initial and
boundary conditions.  However, these  techniques may not be readily extended for
problems  where either the equations have missing parameters or no sufficient
data is available to form a correct initial-value problem.  Such problems are
encountered in various fields of science and engineering and they are handled by
various methods such as PDE-constrained
optimization~\cite{gunzburger2002perspectives}, data
assimilation~\cite{lewis2006dynamic}, optical flow in computer
vision~\cite{fleet2006optical}, and system
identification~\cite{ljung1999system}.  There is a strong need to develop
efficient methods for solving such equations efficiently.

One approach for solving such problems represents the solution using neural
networks.  This approach was introduced by Lagaris et
al.~\cite{lagaris1998artificial}.  At the time, the method did not offer
significant advantages for solving  ODEs and PDEs in  particular due to its
computational cost.  Neural networks were also applied for back-tracking of an
$N$-body system~\cite{quito2001solving}.  Twenty years later, the method was
revived by Raissi et al.~\cite{raissi2019physics} who used modern machine
learning methods and software (such as deep neural networks, automatic
differentiation, and TensorFlow) to carry out its operations. The method has
been popularized by the term physics-informed neural network (PINN).

PINNs cannot match standard numerical methods for solving classical well-posed
problems involving PDEs but they were  positioned as a universal and convenient
tool for solving ill-posed and inverse problems.  However, to the best of our
knowledge, there is no  baseline for comparison and assessment of the
capabilities of these methods.  While the authors admit that conventional
numerical solvers are more efficient for classical
problems~\cite{raissi2019physics}, a direct comparison with such solvers reveals
important drawbacks of the neural network approach.

In PINNs, the cost of evaluating the solution at one point is  proportional to
the number of weights of the neural network, while for conventional grid-based
methods the cost is constant.  Furthermore, as each weight of the neural network
affects the solution at all points, the method is not consistent with the
locality of the original differential problem. This is a  crucial issue as in
PINNs, even for linear problems, the solution is approximated by nonlinear
neural networks.  Consequently, the gradient of the residual with respect to the
weights at one point is a dense vector and the Hessian of the loss function is a
dense matrix.  This makes efficient optimization methods such as Newton's method
unsuitable. In contrast, for linear problems Newton's method would find the
solution exactly after one step and for nonlinear problems would converge
quadratically and typically require few iterations to reach the desired
accuracy.  The supposed convenience of representation and ease of formalism by
PINN comes with a cost.  While PINN evaluates the differential operator exactly
on a set of collocation points, it makes no allowance for physically and
numerically motivated adjustments that speed up convergence and ensure discrete
conservation such as deferred correction, upwind schemes, and discrete
conservation laws as represented in the finite volume method.  Furthermore,
decades of development and analysis that went into conventional solvers enables
understanding, prediction, and control of their convergence and stability
properties. Such information is  not available with neural networks and recent
works are actually aiming to remedy this situation~\cite{mishra2022estimates}
Nevertheless, the PINN framework is actively studied as a tool to  solve
ill-posed and inverse problems~\cite{he2020physics}.

We introduce the Optimizing a Discrete Loss (ODIL) framework that combines
discrete formulations of PDEs with modern machine learning tools to extend
their scope to ill-posed and inverse problems.  The idea of solving the discrete
equations as a minimization problem is known as the
discretize-then-differentiate approach in the context of PDE-constrained
optimization~\cite{gunzburger2002perspectives}, has been formulated for linear
problems as the penalty method~\cite{van2015penalty}, and is related to the
4D-VAR problem in data assimilation~\cite{lewis2006dynamic}.
This procedure has two key aspects. Firstly, the
discretization itself, that defines the accuracy,
stability and consistency of the method.  Secondly, the optimization algorithm
to solve the discrete problem.  Our method uses efficient computational tools
that are  available for both of these aspects.  The discretization properties
are inherited from the classical numerical methods building upon the advances in
this field over multiple decades.  Since the sparsity of the problem is
preserved, the optimization algorithm can use a Hessian and achieve a quadratic
converge rate, which remains out of reach for gradient-based training of neural
networks. Our use of automatic differentiation to compute the Hessian,
makes the implementation as convenient as applying gradient-based methods.

\section{Results}

\subsection{Wave equation. Accuracy and cost of PINN}

\begin{figure}[H]
  \centering
  \includegraphics{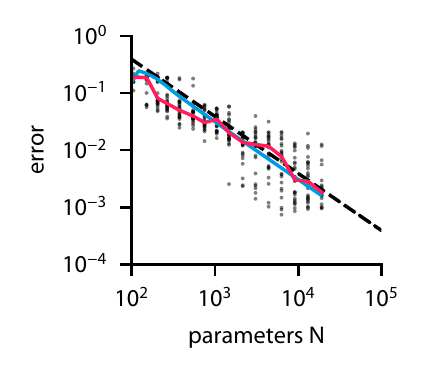}%
  \includegraphics{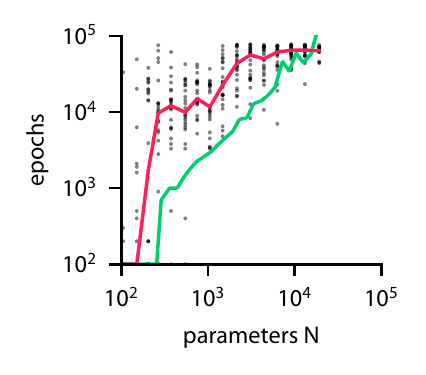}%
  \includegraphics{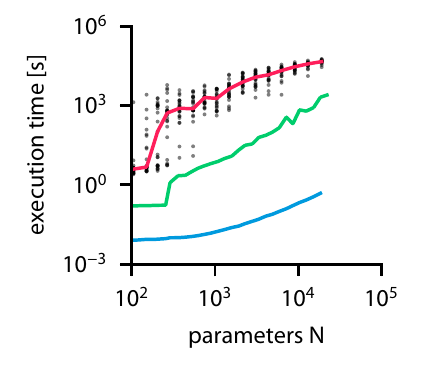}

  \includegraphics{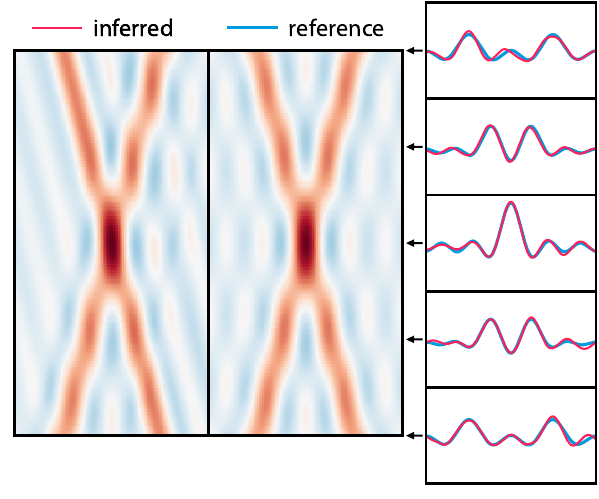}%
  \includegraphics{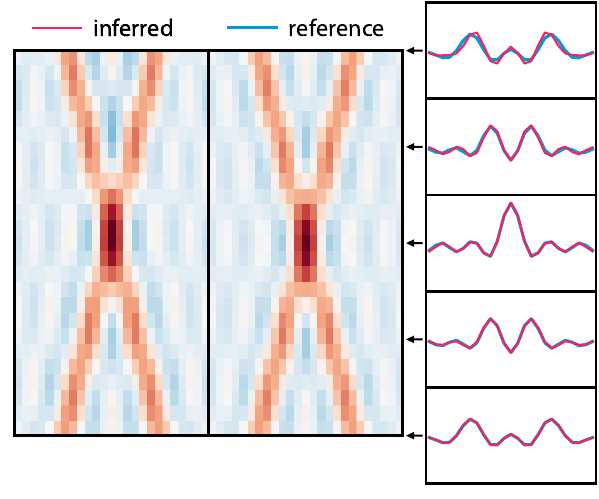}

  \vspace{-7mm}
  \parbox{6cm}{\hspace{9mm}PINN}%
  \parbox{6cm}{\hspace{9mm}ODIL}

  \caption{%
    Wave equation solved using PINN and ODIL.
    Top: $L_2$-error, number of epochs, and execution time of various methods:
    PINN optimized with L-BFGS-B~\key{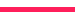},
    ODIL optimized with L-BFGS-B~\key{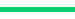},
    ODIL optimized with Newton~\key{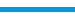},
    and line with slope~$1/N$~\key{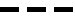}.
    The dots correspond to samples of the initial weights and collocation points
    of PINN. The number of parameters is the number of weights for PINN
    and the number of grid points for ODIL.
    Bottom: Examples of solutions obtained using PINN with two hidden layers of 25
    neurons and ODIL on a $25\times25$ grid.
  }
  \label{f_wave}
\end{figure}

In this section, we apply ODIL to the solution of the  initial-value problem
for the wave equation~$u_{tt}=u_{xx}$. We compare the results of ODIL with those of PINNS in terms of   accuracy and computational cost.
Both methods reduce the problem to a minimization of a loss function. ODIL represents the solution on a uniform grid and constructs the loss
from a second-order finite volume discretization of the equation,
while PINN represents the solution as a fully connected neural network and
evaluates the residuals of the original differential equation exactly on a
finite set of collocation points~\cite{raissi2019physics}.
The initial conditions for $u$ and $u_t$
and Dirichlet conditions on the boundaries are specified from an exact solution
$u(x,t)=\tfrac{1}{10}\sum_{k=1}^5(\cos{(x - t + 0.5)\pi k} + \cos{(x + t +
    0.5)\pi k})$ in a rectangular domain $x\in(-1,1),\;t\in(0,1)$.
For a given number of parameters~$N$, ODIL represents the solution on a
$\sqrt{N}\times\sqrt{N}$ grid, while PINN consists of two equally sized hidden
layers with $\tanh$ activation functions and~$N$ being the total number of
weights.  The number of collocation points for PINN is fixed and amounts to 8192
points inside the domain and 768 points for the initial and boundary conditions.
Figure~\ref{f_wave} shows examples of the solutions obtained using both methods
as well as measurements of their accuracy and cost depending on the number of
parameters.
The measurements for PINN include 25 samples for the random initial weights
and positions of the collocation points.
In the case of ODIL,
solutions from two optimization methods L-BFGS-B and Newton coincide.
The execution time is a product of the number of optimization epochs required
to reach a 150\% of the error obtained after 80000 epochs
and an average execution time over the last 100 epochs.
The measurements are performed on Intel Xeon E5-2690 v3 CPUs.
Both methods demonstrate similar accuracy.
The error of ODIL scales as $1/N$ corresponding to a second-order accurate
discretization.  The error of PINN is significantly scattered, although the
median error also scales as~$1/N$.
The execution times, however, are drastically different.
With $N=20000$ parameters, optimization of PINN using L-BFGS-B takes 15 hours,
optimization of ODIL using L-BFGS-B takes 40 minutes,
and Newton with a direct linear solver~\cite{SciPyNMeth2020}
takes only 0.5 seconds, which is $100\,000$ times faster than PINNs.

\subsection{Lid-driven cavity. Forward problem}

The lid-driven cavity problem is a standard test case~\cite{ghia1982high} for
numerical methods for the steady-state Navier-Stokes equations in two
dimensions:
$u_x + v_y = 0$,
$u u_x + v u_y = -p_x  + 1/\mathrm{Re} (u_{xx} + u_{yy})$,
$u v_x + v v_y = -p_y  + 1/\mathrm{Re} (v_{xx} + v_{yy})$,
where $u(x,y)$ and $v(x,y)$ are the two velocity components and $p(x,y)$ is the
pressure.
The problem is solved in a unit domain with no-slip boundary conditions.
The upper boundary is moving to the right at a unit velocity
while the other boundaries are stagnant.
We use a finite volume discretization on a uniform Cartesian grid with
$128\times128$ cells based on the SIMPLE method~\cite{patankar1983,ferziger2012}
with the Rhie-Chow interpolation~\cite{rhie1983numerical}
to prevent oscillations in the pressure field
and the deferred correction approach
that treats high-order discretization explicitly
and low-order discretization implicitly to obtain an operator with a compact
stencil.

\begin{figure}[tbhp]
  \centering
  \begin{minipage}{4cm}
    \includegraphics[height=20mm]{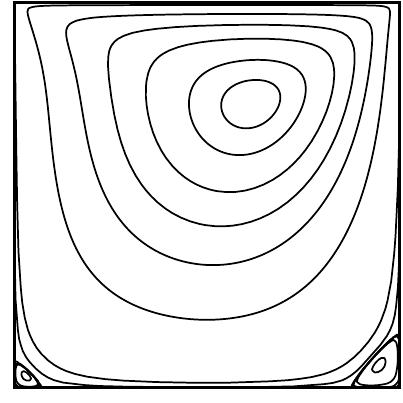}%
    \includegraphics[height=20mm]{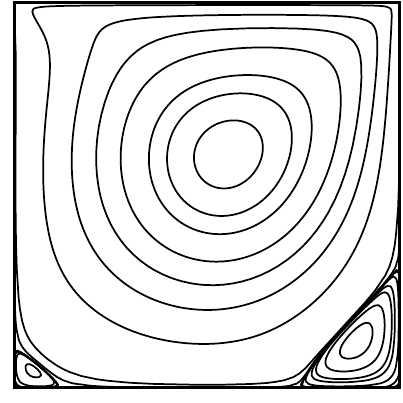}%

    \vspace{-2mm}
    \parbox{20mm}{\quad$\mathrm{Re}=100$}%
    \parbox{20mm}{\quad$\mathrm{Re}=400$}
    \vspace{-1mm}

    \includegraphics[height=20mm]{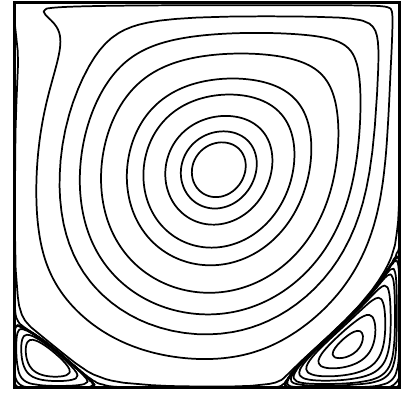}%
    \includegraphics[height=20mm]{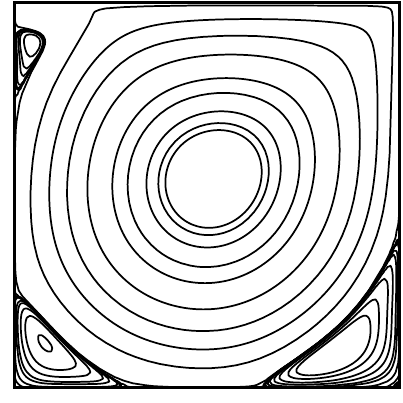}

    \vspace{-2mm}
    \parbox{20mm}{\quad$\mathrm{Re}=1000$}%
    \parbox{20mm}{\quad$\mathrm{Re}=3200$}
  \end{minipage}%
  \begin{minipage}{4cm}
    \hspace*{2mm}\includegraphics{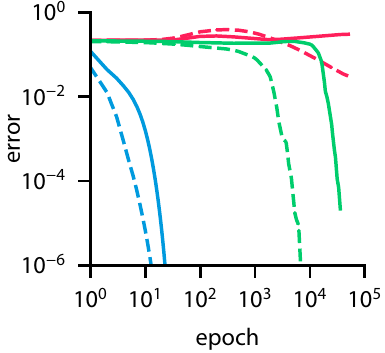}
  \end{minipage}
  \caption{%
    Lid-driven cavity problem solved with ODIL.
    Left: streamlines at various~$\mathrm{Re}$.
    Right: History of the $L_2$-error from various optimization algorithms:
    Adam~\key{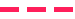}, L-BFGS-B~\key{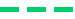}, and Newton~\key{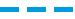} at
    $\mathrm{Re}=100$;
    Adam~\key{110}, L-BFGS-B~\key{120}, and Newton~\key{130} at
    $\mathrm{Re}=400$.
  }
  \label{f_cavity}
\end{figure}

Figure~\ref{f_cavity} shows the streamlines for values of $\mathrm{Re}$
between 100 and 3200, as well as a convergence history of various optimization
algorithms.
The $L_2$ error is computed relative to the solution of the discrete problem.
The number of iterations it takes for L-BFGS-B to reach an error of $10^{-3}$
at $\mathrm{Re}=400$ is in the order of $10000$
which is comparable to the number of iterations typically
needed for training PINNs for flow problems~\cite{raissi2019physics}.
Adam has only reached an error of $0.05$ even after~50\,000
iterations.
Conversely, Newton's method takes about 10 iterations to reach an error of
$10^{-6}$, and the number of iterations does not significantly change
between~$\mathrm{Re}=100$ and $400$.

\subsection{Lid-driven cavity. Flow reconstruction}
\label{s_cavity_reconst}

\begin{figure}[tbhp]
  \centering
  \includegraphics{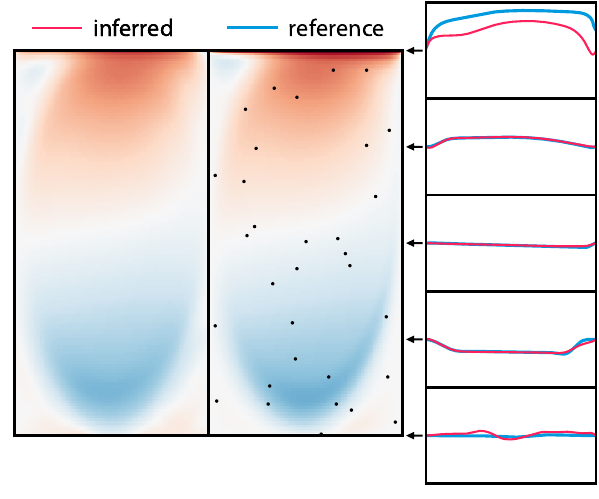}%
  \includegraphics{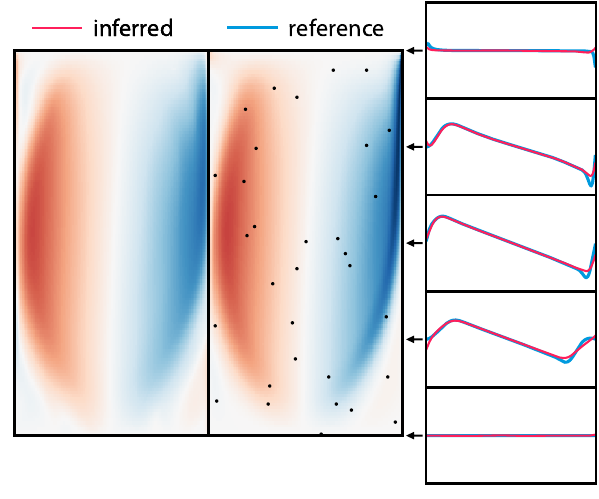}

  \vspace{-7mm}
  \parbox{6cm}{\hspace{9mm}$x$-velocity}%
  \parbox{6cm}{\hspace{9mm}$y$-velocity}

  \caption{%
    Reconstruction of the lid-driven cavity flow at~$\mathrm{Re}=3200$
    from a finite set of points.
    Inferred and reference fields of $x$- and $y$-components of velocity and
    pressure. The reference velocity is imposed in 32 points,
    and the Navier-Stokes equations are imposed everywhere
    with free-slip boundary conditions.
  }
  \label{f_cavity_reconst}
\end{figure}

Here we solve the same equations as in the previous case,
but instead of the no-slip we only impose free-slip boundary conditions
and additionally require that the velocity field equals a reference velocity
field in a finite set of points.
The loss function is a sum of three types of terms:
residuals of equations, terms to impose known velocity values,
and a regularization term
$k^2_\mathrm{reg}\big( \|u_{xx}\|^2+\| u_{yy}\|^2+ \|v_{xx}\|^2+
  \|v_{yy}\|^2\big)$ with $k_\mathrm{reg}=10^{-4}$.
The problem is solved on a $128\times128$ grid and
the reference velocity field is obtained from the forward problem.
Figure~\ref{f_cavity_reconst} shows the results of reconstruction from 32
points. Surprisingly, this small number of points is sufficient to reproduce
features of the flow.

\subsection{Measure of flow complexity}

Based on the flow reconstruction from a finite set of points,
we introduce a measure of flow complexity.
Consider a velocity field $\mb{u}_\mathrm{ref}$ that satisfies the Navier-Stokes
equations.
First we define the quantity
$E(K)=\min\limits_{|X|=K} \| \mb{u}_\mathrm{ref}-\mb{u}(X)\|$
as the best reconstruction error given $K$ points,
where~$\mb{u}(X)$ is the velocity field reconstructed from points $X$
and the minimum is taken over all sets of $K$ points.
Then we define a measure of flow complexity for a given accuracy $\varepsilon>0$
as the minimal number of points required to achieve that accuracy
$K_\mathrm{min}(\varepsilon)=\{K \;|\;E(K)<\varepsilon \}$.

\begin{figure}[H]
  \centering
  \includegraphics{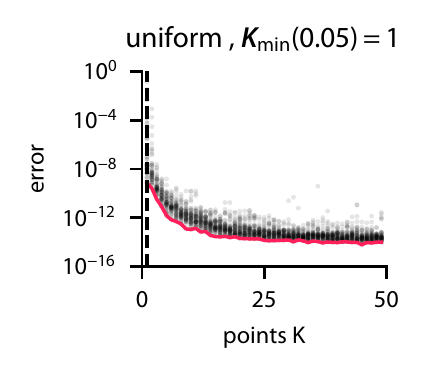}%
  \includegraphics{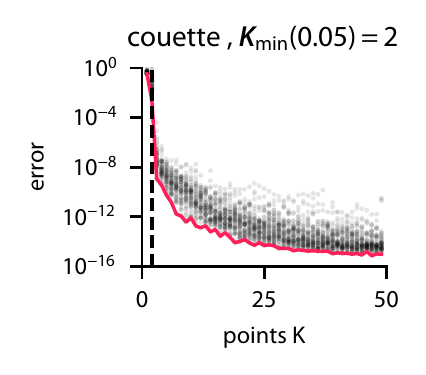}

  \includegraphics{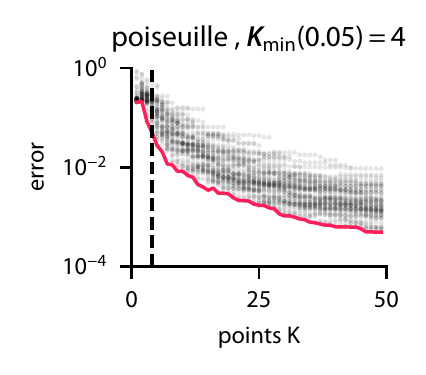}%
  \includegraphics{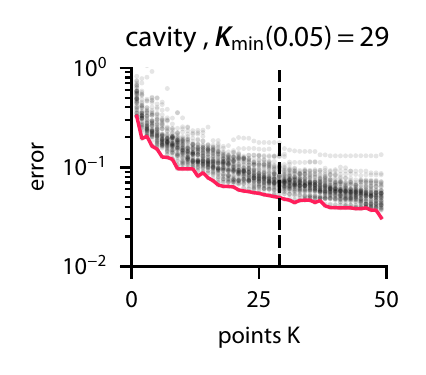}
  \caption{%
    Fluid flow reconstruction error depending on the number
    of reconstruction points
    for various flows: uniform velocity, Couette flow, Poiseuille flow, and the
    flow in a lid-driven cavity.
    The dots correspond to samples of sets of $K$ points
    and the minimal value over them is an estimate of $E(K)$~\key{110}.
  }
  \label{f_measure}
\end{figure}

To illustrate the measure, we consider four types of flow:
uniform velocity, Couette flow (linear profile), Poiseuille flow (parabolic
profile), and the flow in a lid-driven
cavity.  The reference velocity fields in the first three cases are chosen
once at an arbitrary orientation and such that the maximum velocity magnitude is unity.
The problem is solved on a $64\times64$ grid, the velocity at the boundaries is
extrapolated from the inside with second order.
The error is defined through the $L_1$ norm.
The loss function includes the same regularization term as
in~\ref{s_cavity_reconst} with $k_\mathrm{reg}=10^{-3}$.
Figure~\ref{f_measure} shows the values of $E(K)$ and $K_\mathrm{min}(0.05)$.
As expected, the uniform flow takes 1 point to reconstruct
while adding higher order terms to the velocity increases the number of points.
Finally, 29 points are required to reach an error of $0.05$ for the lid-driven
cavity flow.

\subsection{Velocity from tracer}

\begin{figure}[tbhp]
  \centering
  \includegraphics{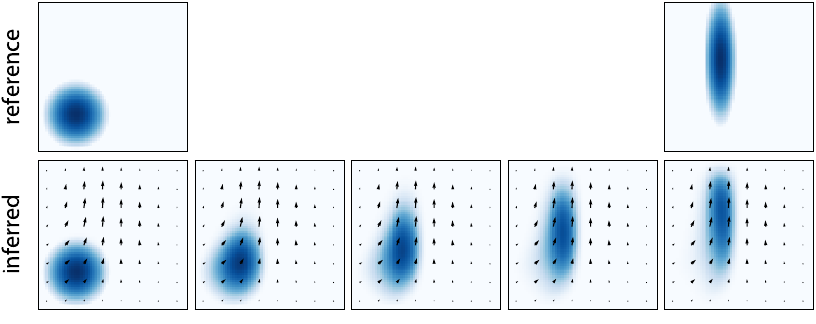}

  \vspace{-2mm}
  \hspace*{2mm}
  \parbox{16mm}{\centering$t=0$}%
  \parbox{16mm}{\centering$0.25$}%
  \parbox{16mm}{\centering$0.5$}%
  \parbox{16mm}{\centering$0.75$}%
  \parbox{16mm}{\centering$1$}%
  \caption{%
    Inferring the velocity field from two snapshots of a tracer.
    The tracer satisfies the advection equation with an unknown steady velocity
    field. The initial and final profiles of the tracer are given.
    Reference tracer field (top) and inferred tracer overlaid by velocity arrows
    (bottom).
  }
  \label{f_tracer}
\end{figure}

We consider the evolution of a tracer field
governed by the advection equation in two dimensions.
The problem is to find a velocity field~$\mb{u}(x, y)$
given that the tracer field~$c(\mb{x},t)$ satisfies the advection equation
$\frac{\partial c}{\partial t} + \mb{u}\cdot\nabla{\alpha} = 0$ in a unit domain,
and the initial~$c(\mb{x},0)=c_0(\mb{x})$
and final~$c(\mb{x},1)=c_1(\mb{x})$ profiles are known.
The discrete problem is solved in space and time on a $64\times64\times64$ grid.
The loss function includes a discretization of the equation with a first-order
upwind scheme,
terms to impose known initial and final profiles of the tracer,
and regularization terms $\|10^{-3}\nabla^2 \mb{u}\|^2$ and $\|\mb{u}_t\|^2$
to prioritize velocity fields that are smooth and stationary.
The results of inference are shown in Figure~\ref{f_tracer}.
The inferred velocity field stretches
the initial tracer representing a circle to match the final profile.

Inferring the velocity field from tracers is an example of
the optical flow problem~\cite{fleet2006optical,mang2018pde} which is important in computer
vision.  Furthermore, reconstructing the velocity field from a concentration
field can assist experimental measurements.

\section{Discussion}
We present a framework (ODIL) for computing the solutions of PDEs by casting their discretization as an optimization problem and applying optimization techniques that are widely available in machine learning software libraries.
The concept of casting the PDE as a loss function has similarities with  similarities with neural network formulations such as PINNs. However the fact that we use the discrete approximation of the equations allows for ODIL  to be  orders of magnitude more
efficient in terms of computational cost and accuracy compared to the PINNs for which complex flow problems ``remain elusive''~\cite{wang2022respecting}.
The present results suggest that ODIL can be a suitable candidate for solving large scale problems
such as weather prediction, complex biological flows, and other scientific and industrial simulations.

\section{Methods}
In the following, we formulate the ODIL framework
for the one-dimensional wave equation $u_{tt}=u_{xx}$
discretized with finite differences.
Instead of solving this equation by marching in time
from known initial conditions, we rewrite the problem as a minimization
of the loss function
\begin{equation*}
  \begin{aligned}
    L(u)
    &=
    \sum\limits_{(i,n)\in \Omega_h}
    \Big(
      \frac{u^{n+1}_i - 2u^n_i + u^{n-1}_i}{\Delta t^2} -
      \frac{u^n_{i+1} - 2u^n_i + u^n_{i-1}}{\Delta x^2}
    \Big)^2
    \\
    &+
    \sum\limits_{(i,n)\in \Gamma_h} \big( u^n_i - g^n_i \big)^2,
  \end{aligned}
  \label{e_waveloss}
\end{equation*}
where the unknown parameter is the solution itself, a discrete
field~$u_i^n$ on a Cartesian grid.
The loss function includes the residual of the discrete equation
in grid points~$\Omega_h$
and imposes boundary and initial conditions~$u=g$ in grid points~$\Gamma_h$.
To solve this problem with a gradient-based method, such as
Adam~\cite{kingma2014adam} or L-BFGS-B~\cite{zhu1997algorithm},
we only need to compute the gradient of the loss
which is also a discrete field~$\pd L/\pd u$.
To apply Newton's method, we assume that the loss function is a sum of
quadratic terms such as
$L(u)=\|\mb{F}[u]\|_2^2+\|\mb{G}[u]\|_2^2$
with discrete operators $\mb{F}[u]$ and $\mb{G}[u]$
and linearize the operators about the current solution $u^s$
to obtain a quadratic loss
\begin{equation*}
  L^s(u^{s+1})=
  \big\|\mb{F}^s+(\pd\mb{F}^s/\pd u) (u^{s+1} - u^s)\big\|_2^2
  +
  \big\|\mb{G}^s+(\pd\mb{G}^s/\pd u) (u^{s+1} - u^s)\big\|_2^2,
\end{equation*}
where $\mb{F}^s$ and $\mb{G}^s$ denote~$\mb{F}[u^s]$ and $\mb{G}[u^s]$.
A minimum of this function provides the solution~$u^{s+1}$ at the next iteration
and satisfies a linear system
\begin{equation*}
  \Big(
  (\pd\mb{F}^s/\pd u)^T (\pd\mb{F}^s/\pd u)
  +
  (\pd\mb{G}^s/\pd u)^T (\pd\mb{G}^s/\pd u)
  \Big)
  (u^{s+1} - u^s)
  +
  (\pd\mb{F}^s/\pd u)^T\mb{F}^s
  +
  (\pd\mb{G}^s/\pd u)^T\mb{G}^s
  =0.
\end{equation*}
Cases with terms other than $\mb{F}$ and $\mb{G}$ are handled similarly.
We further assume that $\mb{F}[u]$ and $\mb{G}[u]$ at each grid point
depend only on the value of $u$ in the neighboring points.
This makes the derivatives $\pd\mb{F}/\pd u$ and $\pd\mb{G}/\pd u$ sparse
matrices.
To implement this procedure, we use automatic
differentiation in TensorFlow~\cite{tensorflow2015whitepaper}
and solve the linear system with either a direct~\cite{SciPyNMeth2020}
or multigrid sparse linear solver~\cite{BeOlSc2022}.

\subsection{Code availability}

The results are obtained using software
available at~\url{https://github.com/cselab/odil}.

\bibliography{mainbib}

\begin{thebibliography}{10}

\bibitem{tensorflow2015whitepaper}
{\sc Abadi, M., Agarwal, A., Barham, P., Brevdo, E., Chen, Z., Citro, C.,
  Corrado, G.~S., Davis, A., Dean, J., Devin, M., Ghemawat, S., Goodfellow, I.,
  Harp, A., Irving, G., Isard, M., Jia, Y., Jozefowicz, R., Kaiser, L., Kudlur,
  M., Levenberg, J., Man\'{e}, D., Monga, R., Moore, S., Murray, D., Olah, C.,
  Schuster, M., Shlens, J., Steiner, B., Sutskever, I., Talwar, K., Tucker, P.,
  Vanhoucke, V., Vasudevan, V., Vi\'{e}gas, F., Vinyals, O., Warden, P.,
  Wattenberg, M., Wicke, M., Yu, Y., and Zheng, X.}
\newblock {TensorFlow}: Large-scale machine learning on heterogeneous systems,
  2015.
\newblock Software available from tensorflow.org.

\bibitem{BeOlSc2022}
{\sc Bell, N., Olson, L.~N., and Schroder, J.}
\newblock {PyAMG}: Algebraic multigrid solvers in python.
\newblock {\em Journal of Open Source Software 7}, 72 (2022), 4142.

\bibitem{ferziger2012}
{\sc Ferziger, J.~H., and Peric, M.}
\newblock {\em Computational methods for fluid dynamics}.
\newblock Springer Science \& Business Media, 2012.

\bibitem{fleet2006optical}
{\sc Fleet, D., and Weiss, Y.}
\newblock Optical flow estimation.
\newblock In {\em Handbook of mathematical models in computer vision}.
  Springer, 2006, pp.~237--257.

\bibitem{ghia1982high}
{\sc Ghia, U., Ghia, K.~N., and Shin, C.}
\newblock {High-Re} solutions for incompressible flow using the {Navier-Stokes}
  equations and a multigrid method.
\newblock {\em Journal of computational physics 48}, 3 (1982), 387--411.

\bibitem{gunzburger2002perspectives}
{\sc Gunzburger, M.~D.}
\newblock {\em Perspectives in flow control and optimization}.
\newblock SIAM, 2002.

\bibitem{he2020physics}
{\sc He, Q., Barajas-Solano, D., Tartakovsky, G., and Tartakovsky, A.~M.}
\newblock Physics-informed neural networks for multiphysics data assimilation
  with application to subsurface transport.
\newblock {\em Advances in Water Resources 141\/} (2020), 103610.

\bibitem{kingma2014adam}
{\sc Kingma, D.~P., and Ba, J.}
\newblock Adam: A method for stochastic optimization.
\newblock {\em arXiv preprint arXiv:1412.6980\/} (2014).

\bibitem{lagaris1998artificial}
{\sc Lagaris, I.~E., Likas, A., and Fotiadis, D.~I.}
\newblock Artificial neural networks for solving ordinary and partial
  differential equations.
\newblock {\em IEEE transactions on neural networks 9}, 5 (1998), 987--1000.

\bibitem{leveque2007finite}
{\sc LeVeque, R.~J.}
\newblock {\em Finite difference methods for ordinary and partial differential
  equations: steady-state and time-dependent problems}.
\newblock SIAM, 2007.

\bibitem{lewis2006dynamic}
{\sc Lewis, J.~M., Lakshmivarahan, S., and Dhall, S.}
\newblock {\em Dynamic data assimilation: a least squares approach}, vol.~13.
\newblock Cambridge University Press, 2006.

\bibitem{ljung1999system}
{\sc Ljung, L.}
\newblock {\em System Identification (2nd Ed.): Theory for the User}.
\newblock Prentice Hall PTR, USA, 1999.

\bibitem{mang2018pde}
{\sc Mang, A., Gholami, A., Davatzikos, C., and Biros, G.}
\newblock Pde-constrained optimization in medical image analysis.
\newblock {\em Optimization and Engineering 19}, 3 (2018), 765--812.

\bibitem{mishra2022estimates}
{\sc Mishra, S., and Molinaro, R.}
\newblock Estimates on the generalization error of physics-informed neural
  networks for approximating a class of inverse problems for pdes.
\newblock {\em IMA Journal of Numerical Analysis 42}, 2 (2022), 981--1022.

\bibitem{patankar1983}
{\sc Patankar, S.~V., and Spalding, D.~B.}
\newblock A calculation procedure for heat, mass and momentum transfer in
  three-dimensional parabolic flows.
\newblock In {\em Numerical Prediction of Flow, Heat Transfer, Turbulence and
  Combustion}. Elsevier, 1983, pp.~54--73.

\bibitem{quito2001solving}
{\sc Quito~Jr, M., Monterola, C., and Saloma, C.}
\newblock Solving n-body problems with neural networks.
\newblock {\em Physical review letters 86}, 21 (2001), 4741.

\bibitem{raissi2019physics}
{\sc Raissi, M., Perdikaris, P., and Karniadakis, G.~E.}
\newblock Physics-informed neural networks: A deep learning framework for
  solving forward and inverse problems involving nonlinear partial differential
  equations.
\newblock {\em Journal of Computational physics 378\/} (2019), 686--707.

\bibitem{rhie1983numerical}
{\sc Rhie, C.~M., and Chow, W.-L.}
\newblock Numerical study of the turbulent flow past an airfoil with trailing
  edge separation.
\newblock {\em AIAA journal 21}, 11 (1983), 1525--1532.

\bibitem{van2015penalty}
{\sc van Leeuwen, T., and Herrmann, F.~J.}
\newblock A penalty method for pde-constrained optimization in inverse
  problems.
\newblock {\em Inverse Problems 32}, 1 (2015), 015007.

\bibitem{SciPyNMeth2020}
{\sc Virtanen, P., Gommers, R., Oliphant, T.~E., Haberland, M., Reddy, T.,
  Cournapeau, D., Burovski, E., Peterson, P., Weckesser, W., Bright, J., {van
  der Walt}, S.~J., Brett, M., Wilson, J., Millman, K.~J., Mayorov, N., Nelson,
  A. R.~J., Jones, E., Kern, R., Larson, E., Carey, C.~J., Polat, {\.I}., Feng,
  Y., Moore, E.~W., {VanderPlas}, J., Laxalde, D., Perktold, J., Cimrman, R.,
  Henriksen, I., Quintero, E.~A., Harris, C.~R., Archibald, A.~M., Ribeiro,
  A.~H., Pedregosa, F., {van Mulbregt}, P., and {SciPy 1.0 Contributors}}.
\newblock {{SciPy} 1.0: Fundamental Algorithms for Scientific Computing in
  Python}.
\newblock {\em Nature Methods 17\/} (2020), 261--272.

\bibitem{wang2022respecting}
{\sc Wang, S., Sankaran, S., and Perdikaris, P.}
\newblock Respecting causality is all you need for training physics-informed
  neural networks.
\newblock {\em arXiv preprint arXiv:2203.07404\/} (2022).

\bibitem{zhu1997algorithm}
{\sc Zhu, C., Byrd, R.~H., Lu, P., and Nocedal, J.}
\newblock Algorithm 778: {L-BFGS-B}: Fortran subroutines for large-scale
  bound-constrained optimization.
\newblock {\em ACM Transactions on mathematical software (TOMS) 23}, 4 (1997),
  550--560.

\bibitem{zienkiewicz2005finite}
{\sc Zienkiewicz, O.~C., Taylor, R.~L., and Zhu, J.~Z.}
\newblock {\em The finite element method: its basis and fundamentals}.
\newblock Elsevier, 2005.

\end{thebibliography}

\end{document}